\documentclass[12pt]{article}
\usepackage{amsfonts,amsmath}
\usepackage{epsfig,graphics}

\renewcommand{\appendix}{%
\renewcommand{\section}{%
\newpage
\thispagestyle{plain}%
\secdef\Appendix\sAppendix}%
\setcounter{section}{0}%
\renewcommand{\thesection}{\Alph{section}}%
}
\newcommand{\Appendix}[2][?]{%
\refstepcounter{section}%
\addcontentsline{toc}{Addendum}%
{\protect\numberline{\appendixname~\thesection}#1}%
{\flushleft\LARGE\bfseries\appendixname\ \thesection\par
\centering#2\par}%
\sectionmark{#1}\vspace{\baselineskip}}

\newcommand{\sAppendix}[1]{%
{\flushright\large\bfseries\appendixname\par
\centering#1\par}%
\vspace{\baselineskip}}
\def\be{\begin{equation}}

\def\bea{\begin{eqnarray}}

\def\eea{\end{eqnarray}}

\normalbaselines \oddsidemargin 0.5cm \evensidemargin 0.5cm
\begin{document}

\pagestyle{empty}

\rightline{ CPHT-RR039.06.08}

\vskip 1cm

\begin{center}

{\Large {\textbf{Multiparameter statistical models from $N^2\times
N^2$ braid matrices: Explicit eigenvalues of transfer matrices
${\bf T}^{(r)}$, spin chains, factorizable scatterings for all
$N$}}}

\vspace{4mm}

{\bf \large B. Abdesselam$^{a,}$\footnote{Email:
boucif@cpht.polytechnique.fr and  boucif@yahoo.fr} and A.
Chakrabarti$^{b,}$\footnote{Email: chakra@cpht.polytechnique.fr}}

\vspace{2mm}

\emph{$^a$ Laboratoire de Physique Th\'eorique, Universit\'e
d'Oran Es-S\'enia, 31100-Oran, Alg\'erie\\
and \\ Laboratoire de Physique Quantique de la Mati\`ere et de
Mod\'elisations Math\'ematiques, Centre Universitaire de Mascara,
29000-Mascara, Alg\'erie}\\
\vspace{2mm}\emph{$^b$ Centre de Physique Th{\'e}orique, Ecole
Polytechnique, 91128 Palaiseau Cedex, France.}
\end{center}

\begin{abstract}
{\small \noindent For a class of multiparameter statistical models
based on $N^2\times N^2$ braid matrices the eigenvalues of the
transfer matrix ${\bf T}^{(r)}$ are obtained explicitly for all
$\left(r,N\right)$. Our formalism yields them as solutions of sets
of linear equations with simple constant coefficients. The role of
zero-sum multiplets constituted in terms of roots of unity is
pointed out and their origin is traced to circular permutations of
the indices in the tensor products of basis states induced by our
class of ${\bf T}^{(r)}$ matrices. The role of free parameters,
increasing as $N^2$ with $N$, is emphasized throughout. Spin chain
Hamiltonians are constructed and studied for all $N$. Inverse
Cayley transforms of Yang-Baxter matrices corresponding to our
braid matrices are obtained for all $N$. They provide potentials
for factorizable $S$-matrices. Main results are summarized and
perspectives are indicated in the concluding remarks.}
\end{abstract}

\pagestyle{plain} \setcounter{page}{1}

\section{Introduction}
\setcounter{equation}{0}

Statistical models "exact" in the sense of Baxter \cite{R1}
satisfy "star-triangle" relations leading to transfer matrices
commuting for different values of the spectral parameter. Crucial
in the study of such models is the spectrum of eigenvalues of
these matrices. But even for the extensively studied 6-vertex and
8-vertex models based on $4\times 4$ braid matrices (the braid
property guaranteeing star-triangle relations), after the first
few simple steps, one has to resort to numerical computations. (We
cannot adequately discuss here the vast associated literature but
refer to texts citing major sources \cite{R2,R3}.) But the basic
reason for such situation is that, starting with $2\times 2$
blocks $T_{ab}$ of the transfer matrix and constructing $2^r\times
2^r$ blocks via coproduct rules for order $r$
$\left(r=1,2,3,\ldots\right)$ one faces increasingly complicated
non-linear systems of equations to be solved in constructing
eigenstates and eigenvalues. Even when a systematic approach is
available, such as the Bethe Ansatz for the 6-vertex case, it only
means that the relevant non-linear equations can be written down
systematically. The task of solving them remains. For the 8-vertex
case one explores analytical properties to extract informations.
(See Ref. 4, for example.) Also the number of free parameters
remains strictly limited for such models, including the multistate
generalization of the 6-vertex one \cite{R2}.

Our class of models exhibits the following properties:

\begin{enumerate}

\item A systematic construction for all dimensions. One starts
from $N^2 \times N^2$ braid matrices leading to $N^r \times N^r$
blocks of ${\bf T}^{(r)}$, the transfer matrix of order $r$ for
$N=2,3,4,5,\ldots$; $r=1,2,3,\ldots$ with no upper limit.

\item The number of free parameters increases with $N$ as $N^2$.
This is a unique feature of our models.

\item For all $\left(N,r\right)$ one solves sets of linear
equations with systematically obtained simple constant
coefficients to construct eigenvalues of ${\bf
T}^{(r)}=\sum_{a=1}^NT_{aa}^{(r)}$. This is a consequence of our
starting point: braid matrices constructed on a basis of "nested
sequence" of projectors \cite{R5}.
\end{enumerate}

The total dimension of the base space on which such equations have
to be solved increases as $N^r$. This will be seen to break up
into subspaces closed under the action of ${\bf T}^{(r)}$, thus
reducing the work considerably. Evidently one cannot continue to
display the results explicitly, as $N^r$ increases. But it is
possible to implement a general and particularly efficient
approach. This will be illustrated for all $N$ and $r=1,2,3,4,5$.
The generalization to $r>5$ will be clearly visible. This is the
central purpose of our paper. Certain other features will be
explored with comments in conclusion.

\section{Braid and transfer matrices}
\setcounter{equation}{0}

The roots of our class of braid matrices and the transfer matrices
they generate are to be found in the nested sequence of projectors
\cite{R5}. In Ref. 5 it was noted that the $4\times 4$ projectors
providing the basis of the 6-vertex and 8-vertex models (i.e. the
braid matrices leading to those models) can be generalized to
higher dimensions -- all higher ones ($N$ odd or even). Braid
matrices for odd dimensions were exhaustively constructed on such
a basis \cite{R6} and related statistical models were studied
\cite{R7}. Then even dimensional ones were presented \cite{R9} and
the corresponding braid matrices and statistical models were
studied \cite{R10}. In previous works the eigenvalues and
eigenfunctions of the transfer matrices were presented mostly for
the lowest values of $N$, namely the $4\times 4$ and $9\times 9$
braid matrices. We present below a general approach for all $N$.

We start by recapitulating the construction of our nested sequence
of projectors, leading to the remarkable properties of our
solutions:

\paragraph{\bf A. Even dimensions:} Let $N=2n$ $\left(n=1,2,3,\ldots\right)$. Define ($\left(ab\right)$ denoting a matrix with
unity on row $a$ and column $b$)
\begin{equation}
P_{ij}^{(\epsilon)}=\frac 12\left\{\left(ii\right)\otimes
\left(jj\right)+\left(\bar{i}\bar{i}\right)\otimes
\left(\bar{j}\bar{j}\right)+\epsilon\left[\left(i\bar{i}\right)\otimes
\left(j\bar{j}\right)+\left(\bar{i}i\right)\otimes
\left(\bar{j}j\right)\right]\right\},
\end{equation}
where $i,j\in\left\{1,\cdots,n\right\}$, $\epsilon=\pm$,
$\bar{i}=2n+1-i$, $\bar{j}=2n+1-j$. Interchanging
$j\rightleftharpoons \bar{j}$ on the right one obtains
$P_{i\bar{j}}^{(\epsilon)}$. One thus obtains a complete basis of
projectors satisfying (with $P_{ij}^{(\epsilon)}=
P_{\bar{i}\bar{j}}^{(\epsilon)}$, $P_{i\bar{j}}^{(\epsilon)}=
P_{\bar{i}j}^{(\epsilon)}$ by definition)
\begin{eqnarray}
&&P_{ab}^{(\epsilon)}P_{cd}^{(\epsilon')}=\delta_{ac}\delta_{bd}\delta_{\epsilon\epsilon'}P_{ab}^{(\epsilon)},\qquad
a,b,c,d\in\left\{1,\ldots,N=2n\right\},\nonumber\\
&&\sum_{\epsilon=\pm}\sum_{i,j=1}^n\left(P_{ij}^{(\epsilon)}+P_{i\bar{j}}^{(\epsilon)}\right)=I_{(2n)^2\times
(2n)^2}.
\end{eqnarray}
Let $m_{ij}^{(\epsilon)}$ be an arbitrary set of parameters
satisfying the crucial constraint
\begin{equation}
m_{ij}^{(\epsilon)}=m_{i\bar{j}}^{(\epsilon)}.
\end{equation}
Define the $N^2\times N^2$ matrix
\begin{equation}
\hat{R}\left(\theta\right)=\sum_{\epsilon}\sum_{i,j}e^{m_{ij}^{(\epsilon)}\theta}\left(P_{ij}^{(\epsilon)}+P_{i\bar{j}}^{(\epsilon)}\right).
\end{equation}
The $\theta$-dependence entering exclusively through the
exponentials as coefficients as above. Such construction
guarantees \cite{R7,R9,R10} the braid property
\begin{equation}
\widehat{R}_{12}\left(\theta\right)\widehat{R}_{23}\left(\theta+\theta'\right)\widehat{R}_{12}\left(\theta'\right)=
\widehat{R}_{23}\left(\theta'\right)\widehat{R}_{12}\left(\theta+\theta'\right)\widehat{R}_{23}\left(\theta\right),
\end{equation}
where, in standard notations, $\widehat{R}_{12}=\widehat{R}\otimes
I$ and $\widehat{R}_{23}=I\otimes \widehat{R}$, $I$ denoting the
$N\times N$ identity matrix.

Define the permutation matrix ${\bf
P}=\sum_{a,b}\left(ab\right)\otimes\left(ba\right)$. Then
\begin{equation}
R\left(\theta\right)={\bf P}\hat{R}\left(\theta\right)
\end{equation}
satisfies the Yang-Baxter equation. The monodromy matrix of order
1 $\left(r=1\right)$ is given by
\begin{equation}
T^{(1)}\left(\theta\right)=R\left(\theta\right).
\end{equation}
The $N\times N$ blocks are (with $a,b\in\left\{1,\ldots,2n
\right\}$)
\begin{equation}
T^{(1)}_{ab}\left(\theta\right)=\frac
12\left(e^{m_{ba}^{(+)}\theta}+e^{m_{ba}^{(-)}\theta}\right)\left(ba\right)+\frac
12\left(e^{m_{ba}^{(+)}\theta}-e^{m_{ba}^{(-)}\theta}\right)\left(\bar{b}\bar{a}\right)\equiv
f_{ba}^{(+)}\left(\theta\right)\left(ba\right)+f_{ba}^{(-)}\left(\theta\right)\left(\bar{b}\bar{a}\right),
\end{equation}
with
\begin{equation}
m_{ba}^{(\epsilon)}=m_{b\bar{a}}^{(\epsilon)}=
m_{\bar{b}a}^{(\epsilon)}=m_{\bar{b}\bar{a}}^{(\epsilon)}.
\end{equation}
Higher orders $\left(r>1\right)$ are obtained via coproduts
\begin{equation}
T^{(r)}_{ab}\left(\theta\right)=\sum_{c_1,\ldots,c_{r-1}}T^{(1)}_{ac_1}\left(\theta\right)\otimes
T^{(1)}_{c_1c_2}\left(\theta\right)\otimes\ldots\otimes
T^{(1)}_{c_{r-2}c_{r-1}}\left(\theta\right)\otimes
T^{(1)}_{c_{r-1}b}\left(\theta\right).
\end{equation}
For some essential purposes it is worthwhile express (2.8) as
(with $i,j\in\left\{1,\ldots,n\right\}$)
\begin{eqnarray}
&&T^{(1)}_{ij}\left(\theta\right)=\frac 12\sum_{\epsilon}
e^{m_{ji}^{(\epsilon)}\theta}\left(\left(ji\right)+\epsilon\left(\bar{j}\bar{i}\right)\right),\nonumber\\
&&T^{(1)}_{\bar{i}\bar{j}}\left(\theta\right)=\frac
12\sum_{\epsilon}
e^{m_{ji}^{(\epsilon)}\theta} \left(\left(\bar{j}\bar{i}\right)+\epsilon\left(ji\right)\right),\nonumber\\
&&T^{(1)}_{i\bar{j}}\left(\theta\right)=\frac 12\sum_{\epsilon}
e^{m_{ji}^{(\epsilon)}\theta}\left(\left(\bar{j}i\right)+\epsilon\left(j\bar{i}\right)\right),\nonumber\\
&&T^{(1)}_{\bar{i}j}\left(\theta\right)=\frac 12\sum_{\epsilon}
e^{m_{ji}^{(\epsilon)}\theta}\left(\left(j\bar{i}\right)+\epsilon\left(\bar{j}i\right)\right).
\end{eqnarray}
The transfer matrix is the trace on $\left(a,b\right)$:
\begin{equation}
{\bf T}^{(r)}\left(\theta\right)=\sum_a
T^{(r)}_{aa}\left(\theta\right).
\end{equation}
The foregoing construction assures the crucial commutativity
\begin{equation}
\left[{\bf T}^{(r)}\left(\theta\right),{\bf
T}^{(r)}\left(\theta'\right)\right]=0.
\end{equation}
This implies that the eigenstates of ${\bf
T}^{(r)}\left(\theta\right)$ are $\theta$-independent (sum of
basis vectors with $\theta$-independent, constant relative
coefficients).

\paragraph{\bf B. Odd dimensions:} From the extensive previous studies
\cite{R6,R7} we select the essential features. For $N=2n-1$
($n=2,3,\ldots$)
\begin{equation}
\bar{n}=N-n+1=n
\end{equation}
with the same definition as in (2.1). This is the crucial new
feature. For $\left(i,j\right)\neq n$, one has the same
$P_{ij}^{(\epsilon)}$, $P_{i\bar{j}}^{(\epsilon)}$ as before. But
now
\begin{eqnarray}
&&P_{in}^{(\epsilon)}=\frac
12\left\{\left(ii\right)+\left(\bar{i}\bar{i}\right)+\epsilon\left[\left(i\bar{i}\right)+\left(\bar{i}i\right)\right]\right\}\otimes
\left(nn\right),\nonumber\\
&&P_{ni}^{(\epsilon)}=\frac 12
\left(nn\right)\otimes \left\{\left(ii\right)+\left(\bar{i}\bar{i}\right)+\epsilon\left[\left(i\bar{i}\right)+\left(\bar{i}i\right)
\right]\right\},\nonumber\\
&&P_{nn}^{(\epsilon)}=\left(nn\right)\otimes \left(nn\right).
\end{eqnarray}
Normalizing the coefficient of $P_{nn}^{(\epsilon)}$ to unity
$\hat{R}\left(\theta\right)$ of the (2.4) now becomes
\begin{equation}
\hat{R}\left(\theta\right)=P_{nn}^{(\epsilon)} +\sum_{i,\epsilon}
\left(e^{m_{ni}^{(\epsilon)}\theta}P_{ni}^{(\epsilon)}+e^{m_{in}^{(\epsilon)}\theta}P_{in}^{(\epsilon)}\right)
+\sum_{i,j\epsilon}e^{m_{ij}^{(\epsilon)}\theta}\left(P_{ij}^{(\epsilon)}+P_{i\bar{j}}^{(\epsilon)}\right)
\end{equation}
($i,j\in\left\{1,\cdots,n-1\right\}$,
$\bar{i},\bar{j}\in\left\{2n-1,\cdots,n+1\right\}$,
$\epsilon=\pm$). There are evident parallel modifications
concerning $T^{(r)}_{ab}\left(\theta\right)$. Ref. 7 contains
detailed discussions and results for $N=3$, $r=1,2,3,4$. We will
not repeat them here, but will reconsider them later in the
context of the general methods to be implemented below.

\paragraph{\bf C. Free parameters:} A unique feature of our
constructions is the number of free parameters
$m_{ij}^{(\epsilon)}$, increasing as $N^2$ with $N$. For our
choice of normalizations (apart from a possible altered choice of
an overall factor, irrelevant for (2.5)) the exact numbers are
\begin{eqnarray}
&&\frac 12N^2=2n^2\qquad\hbox{for}\qquad N=2n,\nonumber\\
&&\frac
12\left(N+3\right)\left(N-1\right)=2\left(n^2-1\right)\qquad\hbox{for}\qquad
N=2n-1.
\end{eqnarray}
This is to be contrasted with the multistate generalization of the
6-vertex model where the parametrization remain fixed at the
6-vertex level (Ref. 2 and sources cited here).

One of our principal aims is to display the roles of our free
parameters concerning the basic features of our models.

\section{Explicit eigenvalues of Transfer Matrices for all dimensions $N$ and orders $r$}
\setcounter{equation}{0}

We present below a unified approach for all $N$  and illustrate it
in some detail for $r=1,2,3,4,5$. For $r>5$ the extension will be
evident.

\paragraph{\bf A. Even dimensions:} In sec. 5 of Re. 9 it was pointed
out that the antidiagonal matrix
\begin{equation}
K=\sum_{a=1}^{2n}\left(a\bar{a}\right)=\sum_{i=1}^{n}\left(\left(i\bar{i}\right)+\left(\bar{i}i\right)\right)
\end{equation}
relates the blocks of $T^{(1)}$ as follows:
\begin{equation}
KT^{(1)}_{ab}=T^{(1)}_{a\bar{b}},\qquad
T^{(1)}_{ab}K=T^{(1)}_{\bar{a}b},\qquad
KT^{(1)}_{ab}K=T^{(1)}_{\bar{a}\bar{b}}.
\end{equation}
For $N=2$ this led to an iterative construction of eigenstates and
eigenvalues of the transfer matrix ${\bf T}^{(r)}$ for increasing
$r$. This works efficiently only for $N=2$ (the $4\times 4$
$\hat{R}\left(\theta\right)$-matrix).

For our present approach the crucial ingredient shows up
immediately in the structure (arising from those of our projectors
and the coproducts rules)
\begin{eqnarray}
&&{\bf T}^{(r)}=\displaystyle\frac 1{2^r}\sum_{a_1,a_2,\ldots,a_r}
\sum_{\epsilon_{21},\epsilon_{32},\ldots,\epsilon_{1r}}e^{(m_{a_2a_1}^{(\epsilon_{21})}+m_{a_3a_2}^{(\epsilon_{32})}
+\cdots+m_{a_1a_r}^{(\epsilon_{1r})})\theta}\nonumber\\
&&\phantom{{\bf T}^{(r)}=}
\left[\left(a_2a_1\right)+\epsilon_{21}\left(\bar{a_2}\bar{a_1}\right)\right]\otimes
\left[\left(a_3a_2\right)+\epsilon_{32}\left(\bar{a_3}\bar{a_2}\right)\right]
\otimes\cdots\otimes
\left[\left(a_1a_r\right)+\epsilon_{1r}\left(\bar{a_1}\bar{a_r}\right)\right],
\end{eqnarray}
where each $\epsilon$ is independently $\left(\pm\right)$. From
(3.3), one obtains immediately,
\begin{eqnarray}
&&\left({\bf T}^{(r)}\right)^{-1}=\frac
1{2^r}\sum_{a_1,a_2,\ldots,a_r}
\sum_{\epsilon_{21},\epsilon_{32},\ldots,\epsilon_{1r}}e^{-(m_{a_2a_1}^{(\epsilon_{21})}+m_{a_3a_2}^{(\epsilon_{32})}
+\cdots+m_{a_1a_r}^{(\epsilon_{1r})})\theta}\\
&&\phantom{\left({\bf T}^{(r)}\right)^{-1}=}
\left[\left(a_1a_2\right)+\epsilon_{21}\left(\bar{a_1}\bar{a_2}\right)\right]\otimes
\left[\left(a_2a_3\right)+\epsilon_{32}\left(\bar{a_2}\bar{a_3}\right)\right]
\otimes\cdots\otimes
\left[\left(a_ra_1\right)+\epsilon_{1r}\left(\bar{a_r}\bar{a_1}\right)\right].\nonumber
\end{eqnarray}
Consider the basis state
\begin{equation}
\left|b_1\right\rangle\otimes\left|b_2\right\rangle\otimes\ldots\otimes
\left|b_r\right\rangle\equiv\left|b_1b_2\ldots b_r\right\rangle
\end{equation}

\paragraph{\bf (1):} If any one (even a single one) of the indices
$b_i\neq a_{i+1}$ or $\overline{a_{i+1}}$ (cyclic) then the
coefficient of
$e^{\left(m_{a_2a_1}^{(\epsilon_{21})}+\cdots+m_{a_1a_r}^{(\epsilon_{1r})}\right)\theta}$
namely
$\left[\left(a_2a_1\right)+\epsilon_{21}\left(\bar{a_2}\bar{a_1}\right)\right]\otimes
\cdots\otimes
\left[\left(a_1a_r\right)+\epsilon_{1r}\left(\bar{a_1}\bar{a_r}\right)\right]$
will annihilate it. On the other hand (2.9) plays a crucial role.
Gathering together all the terms with
$e^{\left(m_{b_2b_1}^{(\epsilon_{21})}+m_{b_3b_2}^{(\epsilon_{32})}
+\cdots+m_{b_rb_{r-1}}^{(\epsilon_{r,r-1})}+m_{b_1b_r}^{(\epsilon_{1r})}\right)\theta}$
as coefficient one obtains the action of ${\bf T}^{(r)}$ on (3.5).
All essential features for the general case can be read off the
first few simples examples:
\begin{eqnarray}
&&2^2{\bf T}^{(2)}\left|b_1b_2\right\rangle=
\sum_{\epsilon_{21},\epsilon_{12}}e^{\left(m_{b_2b_1}^{(\epsilon_{21})}+m_{b_1b_2}^{(\epsilon_{12})}\right)\theta}
\left(1+\epsilon_{21}\epsilon_{12}\right)\left[\left|b_2b_1\right\rangle+\epsilon_{21}\left|\bar{b_2}\bar{b_1}\right\rangle\right],\\
&&2^3{\bf T}^{(3)}\left|b_1b_2b_3\right\rangle=
\sum_{\epsilon_{21},\epsilon_{32},\epsilon_{13}}e^{\left(m_{b_2b_1}^{(\epsilon_{21})}+m_{b_3b_2}^{(\epsilon_{32})}+
m_{b_1b_3}^{(\epsilon_{13})}\right)\theta}
\left(1+\epsilon_{21}\epsilon_{32}\epsilon_{13}\right)\nonumber\\
&&\left[\left|b_2b_3b_1\right\rangle+\epsilon_{21}\left|\bar{b_2}b_3\bar{b_1}\right\rangle+
\epsilon_{32}\left|\bar{b_2}\bar{b_3}b_1\right\rangle+\epsilon_{13}\left|b_2\bar{b_3}\bar{b_1}\right\rangle\right],\\
&&2^4{\bf T}^{(4)}\left|b_1b_2b_3b_4\right\rangle=
\sum_{\epsilon_{21},\epsilon_{32},\epsilon_{43},\epsilon_{14}}e^{\left(m_{b_2b_1}^{(\epsilon_{21})}+m_{b_3b_2}^{(\epsilon_{32})}+
m_{b_4b_3}^{(\epsilon_{43})}+m_{b_1b_4}^{(\epsilon_{14})}\right)\theta}\left(1+\epsilon_{21}\epsilon_{32}\epsilon_{43}\epsilon_{14}\right)\nonumber\\
&&\left[\left|b_2b_3b_4b_1\right\rangle+\epsilon_{21}\left|\bar{b_2}b_3b_4\bar{b_1}\right\rangle+
\epsilon_{32}\left|\bar{b_2}\bar{b_3}b_4b_1\right\rangle+\epsilon_{43}\left|b_2\bar{b_3}\bar{b_4}b_1\right\rangle+
\right.\nonumber\\
&&\left.\epsilon_{14}\left|b_2b_3\bar{b_4}\bar{b_1}\right\rangle+\epsilon_{21}\epsilon_{32}\left|b_2\bar{b_3}b_4\bar{b_1}\right\rangle+
\epsilon_{21}\epsilon_{43}\left|\bar{b_2}\bar{b_3}\bar{b_4}\bar{b_1}\right\rangle+\epsilon_{21}\epsilon_{14}
\left|\bar{b_2}b_3\bar{b_4}b_1\right\rangle\right]
\end{eqnarray}
Note that the bars correspond to the indies of the $\epsilon$'s.
One has thus $\epsilon_{ij}\left|\bar{b_i}\bar{b_j}
b_kb_l\right\rangle$, $\epsilon_{ij}\epsilon_{ki}
\left|\bar{\bar{b_i}}\bar{b_j}
\bar{b_k}b_l\right\rangle=\epsilon_{ij}\epsilon_{ki}
\left|b_i\bar{b_j} \bar{b_k}b_l\right\rangle$ since
$\bar{\bar{b_i}}=b_i$. Consequently only even number of additional
bars appear on the right. Starting with even or odd number of bars
on the left, lead to "even " and "odd" closed subspaces.
\begin{eqnarray}
&&2^5{\bf T}^{(5)}\left|b_1b_2b_3b_4b_5\right\rangle=
\sum_{\epsilon_{21},\epsilon_{32},\epsilon_{43},\epsilon_{14}}e^{\left(m_{b_2b_1}^{(\epsilon_{21})}+m_{b_3b_2}^{(\epsilon_{32})}+
m_{b_4b_3}^{(\epsilon_{43})}+m_{b_5ba_4}^{(\epsilon_{54})}+m_{b_1b_5}^{(\epsilon_{15})}\right)\theta}
\left(1+\epsilon_{21}\epsilon_{32}\epsilon_{43}\epsilon_{54}\epsilon_{15}\right)\nonumber\\
&&\left[\left|b_2b_3b_4b_5b_1\right\rangle+\epsilon_{21}\left|\bar{b_2}b_3b_4b_5\bar{b_1}\right\rangle+
\epsilon_{32}\left|\bar{b_2}\bar{b_3}b_4b_4b_1\right\rangle+\epsilon_{43}\left|b_2\bar{b_3}\bar{b_4}b_5b_1\right\rangle+
\right.\nonumber\\
&&\left.\epsilon_{54}\left|b_2b_3\bar{b_4}\bar{b_5}b_1\right\rangle+\epsilon_{15}\left|b_2b_3b_4\bar{b_5}\bar{b_1}\right\rangle+
\epsilon_{21}\epsilon_{32}\left|b_2\bar{b_3}b_4b_5\bar{b_1}\right\rangle+
\epsilon_{21}\epsilon_{43}\left|\bar{b_2}\bar{b_3}\bar{b_4}b_5\bar{b_1}\right\rangle+\right.\nonumber\\
&&\left.\epsilon_{21}\epsilon_{54}
\left|\bar{b_2}b_3\bar{b_4}\bar{b_5}\bar{b_1}\right\rangle+\epsilon_{21}\epsilon_{15}
\left|\bar{b_2}b_3b_4\bar{b_5}b_1\right\rangle+\epsilon_{32}\epsilon_{43}
\left|\bar{b_2}b_3\bar{b_4}b_5b_1\right\rangle+\epsilon_{32}\epsilon_{54}
\left|\bar{b_2}\bar{b_3}\bar{b_4}\bar{b_5}b_1\right\rangle+\right.\nonumber\\
&&\left.\epsilon_{32}\epsilon_{15}
\left|\bar{b_2}\bar{b_3}b_4\bar{b_5}\bar{b_1}\right\rangle+\epsilon_{43}\epsilon_{54}
\left|b_2\bar{b_3}b_4\bar{b_5}b_1\right\rangle+\epsilon_{43}\epsilon_{15}
\left|b_2\bar{b_3}\bar{b_4}\bar{b_5}\bar{b_1}\right\rangle+\epsilon_{54}\epsilon_{15}
\left|b_2b_3\bar{b_4}b_5\bar{b_1}\right\rangle\right]
\end{eqnarray}

\paragraph{\bf (2):} It is important to note that the preceding
features are independent of $N$ -- for a given $r$ they are valid
for all $N$. For $r>N$ all the indices $\left(b_1,\ldots,b_r
\right)$ cannot be distinct. But even for $r<N$ some or all of
them can coincide. This general validity permits a unified
treatment for all $N$.

\paragraph{\bf (3):} Note also the crucial cyclic permutation of
the indices under the action of ${\bf T}^{(r)}$:
$\left(b_1,b_2,b_3,\ldots,b_{r-1},b_r \right)\longrightarrow
\left(b_2,b_3,b_4,\ldots,b_{r},b_1 \right)$. This is a consequence
of the trace condition in (2.12) incorporated in (3.3). It will be
seen later how this leads to an essential role of the $r$-th roots
of unity $e^{{\bf i}\frac{k\cdot 2\pi}r}$ ($k=0,1,\ldots,r-1$) in
our construction of eigenstates and also of eigenvalues as factors
accompanying exponentials of the type appearing in (3.6-9).

\vskip 3mm

We continue the study of the particular cases $r=1,2,3,4,5$ each
one for any $N$. General features will be better understood
afterwards.

\paragraph{$\bullet$ $r=1$:}
\begin{equation}
{\bf T}^{(1)}= \frac 12
\sum_{a,\epsilon}e^{m_{aa}^{(\epsilon)}\theta}
\left[\left(aa\right)+\epsilon\left(\bar{a}\bar{a}\right)\right]
=\frac 12
\sum_{i,\epsilon}\left(1+\epsilon\right)e^{m_{ii}^{(\epsilon)}\theta}
\left[\left(ii\right)+\epsilon\left(\bar{i}\bar{i}\right)\right],
\end{equation}
where $a\in\left\{1,\ldots,2n\right\}$, $i\in\left\{1,\ldots,
n\right\}$. This is already in diagonal form. Eigenstates and
values are obtained trivially.

\paragraph{$\bullet$ $r=2$:} We start with (3.6). Separating "even" and "odd" spaces,
define (for $i\neq j$ when $\left(i,j,\bar{i},\bar{j}\right)$ are
all distinct)
\begin{eqnarray}
&&{\bf
T}^{(2)}\left(c_1\left|ij\right\rangle+c_2\left|\bar{i}\bar{j}\right\rangle
+c_3\left|ji\right\rangle+c_4\left|\bar{j}\bar{i}\right\rangle\right)=
\upsilon_e\left(c_1\left|ij\right\rangle+c_2\left|\bar{i}\bar{j}\right\rangle
+c_3\left|ji\right\rangle+c_4\left|\bar{j}\bar{i}\right\rangle\right),\\
&&{\bf
T}^{(2)}\left(d_1\left|i\bar{j}\right\rangle+d_2\left|\bar{i}j\right\rangle
+d_3\left|j\bar{i}\right\rangle+d_4\left|\bar{j}i\right\rangle\right)=
\upsilon_o\left(d_1\left|i\bar{j}\right\rangle+d_2\left|\bar{i}j\right\rangle
+d_3\left|j\bar{i}\right\rangle+d_4\left|\bar{j}i\right\rangle\right).\qquad
\phantom{x}
\end{eqnarray}
The solutions are, for
$\left(c_1,c_2,c_3,c_4\right)=\left(1,1,1,1\right),
\left(1,1,-1,-1\right),\left(1,-1,1,-1\right),\left(1,-1,-1,1\right)$
respectively
\begin{equation}
\upsilon_e=e^{\left(m_{ij}^{(+)}+m_{ji}^{(+)}\right)\theta},\,
-e^{\left(m_{ij}^{(+)}+m_{ji}^{(+)}\right)\theta},\,
e^{\left(m_{ij}^{(-)}+m_{ji}^{(-)}\right)\theta},\,
-e^{\left(m_{ij}^{(-)}+m_{ji}^{(-)}\right)\theta}
\end{equation}
and for $\left(d_1,d_2,d_3,d_4\right)=\left(1,1,1,1\right),\,
\left(1,1,-1,-1\right),\,\left(1,-1,1,-1\right),\,\left(1,-1,-1,1\right)$
respectively
\begin{equation}
\upsilon_o=e^{\left(m_{ij}^{(+)}+m_{ji}^{(+)}\right)\theta},\,
-e^{\left(m_{ij}^{(+)}+m_{ji}^{(+)}\right)\theta},\,
-e^{\left(m_{ij}^{(-)}+m_{ji}^{(-)}\right)\theta},\,
e^{\left(m_{ij}^{(-)}+m_{ji}^{(-)}\right)\theta}
\end{equation}
For $i=j$, with $\epsilon=\pm$
\begin{eqnarray}
&&{\bf
T}^{(2)}\left(\left|ii\right\rangle+\epsilon\left|\bar{i}\bar{i}\right\rangle
\right)=e^{2m_{ii}^{(\epsilon)}\theta}\left(\left|ii\right\rangle+\epsilon\left|\bar{i}\bar{i}\right\rangle
\right),\nonumber\\
&&{\bf
T}^{(2)}\left(\left|i\bar{i}\right\rangle+\epsilon\left|\bar{i}i\right\rangle
\right)=\epsilon
e^{2m_{ii}^{(\epsilon)}\theta}\left(\left|i\bar{i}\right\rangle+\epsilon\left|\bar{i}i\right\rangle
\right).\end{eqnarray} Summing over all $i$, consistently with the
general rule, one obtains
\begin{equation}
\hbox{Tr}\left({\bf T}^{(2)}\right) =
2\sum_{i=1}^ne^{2m_{ii}^{(+)}\theta}
\end{equation}
contributions coming only for $i=j$. For $r=2$, the $r$-th roots
of unity are $\pm 1$ and their role above is not visible very
distinctly due to the presence of $\pm 1$ also from $\epsilon$.
The role of the roots of unity will be more evident from $r=3$
onwards. The simple exercise above gives all eigenvalues and
eigenstates for ${\bf T}^{(2)}$ for all $N(=2n)$. For $N=2$,
$r=2$, writing out the $16\times 16$ ${\bf T}^{(2)}$ in full and
effectively diagonalizing it we have confirmed the results above.

\paragraph{$\bullet$ $r=3$:} Now the starting point is (3.7).
Define, for the even subset,
\begin{eqnarray}
&&V_{\lambda}=\left(\left|i_1i_2i_3\right\rangle+\epsilon_{12}\left|\bar{i_1}\bar{i_2}i_3\right\rangle
+\epsilon_{23}\left|i_1\bar{i_2}\bar{i_3}\right\rangle+\epsilon_{31}\left|\bar{i_1}i_2\bar{i_3}\right\rangle\right)
+\nonumber\\
&&\phantom{V_{\lambda}=}\lambda\left(\left|i_3i_1i_2\right\rangle+\epsilon_{12}\left|i_3\bar{i_1}\bar{i_2}\right\rangle
+\epsilon_{23}\left|\bar{i_3}i_1\bar{i_2}\right\rangle+\epsilon_{31}\left|\bar{i_3}\bar{i_1}i_2\right\rangle\right)+\nonumber\\
&&\phantom{V_{\lambda}=}\lambda^2\left(\left|i_2i_3i_1\right\rangle+\epsilon_{12}\left|\bar{i_2}i_3\bar{i_1}\right\rangle
+\epsilon_{23}\left|\bar{i_2}\bar{i_3}i_1\right\rangle+\epsilon_{31}\left|i_2\bar{i_3}\bar{i_1}\right\rangle\right)
\end{eqnarray}
where
\begin{equation}
\lambda=\left(1,e^{{\bf i}\frac{2\pi}3} ,e^{{\bf
i}\frac{4\pi}3}\right)
\end{equation}
and hence $\lambda^3=1$, and also $1+e^{{\bf
i}\frac{2\pi}3}+e^{{\bf i}\frac{4\pi}3}=0$. Corresponding to
$\left(i_1,i_2,i_3\right)\rightleftharpoons
\left(\bar{i_1},\bar{i_2},\bar{i_3}\right)$, define for the
odd-subset,
\begin{eqnarray}
&&\overline{V_{\lambda}}=\left(\left|\bar{i_1}\bar{i_2}\bar{i_3}\right\rangle+\epsilon_{12}\left|i_1i_2\bar{i_3}\right\rangle
+\epsilon_{23}\left|\bar{i_1}i_2i_3\right\rangle+\epsilon_{31}\left|i_1\bar{i_2}i_3\right\rangle\right)
+\nonumber\\
&&\phantom{V_{\lambda}=}\lambda\left(\left|\bar{i_3}\bar{i_1}\bar{i_2}\right\rangle+\epsilon_{12}\left|\bar{i_3}i_1i_2\right\rangle
+\epsilon_{23}\left|i_3\bar{i_1}i_2\right\rangle+\epsilon_{31}\left|i_3i_1\bar{i_2}\right\rangle\right)+\nonumber\\
&&\phantom{V_{\lambda}=}\lambda^2\left(\left|\bar{i_2}\bar{i_3}\bar{i_1}\right\rangle+\epsilon_{12}\left|i_2\bar{i_3}i_1\right\rangle
+\epsilon_{23}\left|i_2i_3\bar{i_1}\right\rangle+\epsilon_{31}\left|\bar{i_2}i_3i_1\right\rangle\right).
\end{eqnarray}
The 8 values of the set $\left(\epsilon_{12},\epsilon_{23},
\epsilon_{31}\right)$ are restricted (see (3.7)) by $\left(1+
\epsilon_{12} \epsilon_{23} \epsilon_{31}\right)=2$ ($\neq 0$).
The remaining 4 possibilities along with 3 for $\lambda$ (see
(3.18)) yields $4\times 3=12$ solutions for $V_\lambda$ and
$\overline{V_\lambda}$ each for distinct indices $i_1\neq i_2\neq
i_3\neq i_1$. The eigenvalues are both $V_\lambda$ and
$\overline{V_\lambda}$,
\begin{equation}
e^{\left(m_{i_1i_2}^{(\epsilon_{12})}+m_{i_2i_3}^{(\epsilon_{23})}+m_{i_3i_1}^{(\epsilon_{31})}\right)\theta}\left(1,e^{{\bf
i}\frac{2\pi}3} ,e^{{\bf i}\frac{4\pi}3}\right).
\end{equation}
For two of the three indices equal ($i_1=i_2\neq i_3$,$\ldots$)
the preceding results can be carried over without crucial change.
But for $i_1=i_2=i_3\equiv i$, say, there is a special situation
due to the fact
\begin{equation}
1+\lambda+\lambda^2=\left(3,0,0\right).
\end{equation}
for $\lambda=\left(1,e^{{\bf i}\frac{2\pi}3} ,e^{{\bf
i}\frac{4\pi}3}\right)$ respectively. Now \begin{eqnarray}
&&V_{\lambda}=\left(1+\lambda+\lambda^2\right)\left|iii\right\rangle
+\left(\epsilon_{12}+\epsilon_{31}\lambda+\epsilon_{23}\lambda^2\right)\left|\bar{i}\bar{i}i\right\rangle+\nonumber\\
&&\phantom{V_{\lambda}=}\left(\epsilon_{23}+\epsilon_{12}\lambda+\epsilon_{31}\lambda^2\right)\left|i\bar{i}\bar{i}\right\rangle+
\left(\epsilon_{31}+\epsilon_{23}\lambda+\epsilon_{12}\lambda^2\right)\left|\bar{i}i\bar{i}\right\rangle.\end{eqnarray}
The eigenfunctions and eigenvalues are (for mutually orthogonal
basis states)
\begin{eqnarray}
&&{\bf T}^{(3)} \left(\left|iii\right\rangle
+\left|\bar{i}\bar{i}i\right\rangle+\left|i\bar{i}\bar{i}\right\rangle+\left|\bar{i}i\bar{i}\right\rangle\right)
=e^{3m_{ii}^{(+)}\theta}\left(\left|iii\right\rangle
+\left|\bar{i}\bar{i}i\right\rangle+\left|i\bar{i}\bar{i}\right\rangle+\left|\bar{i}i\bar{i}\right\rangle\right)\\
&&{\bf T}^{(3)} \left(3\left|iii\right\rangle
-\left|\bar{i}\bar{i}i\right\rangle-\left|i\bar{i}\bar{i}\right\rangle-\left|\bar{i}i\bar{i}\right\rangle\right)
=e^{\left(m_{ii}^{(+)}+2m_{ii}^{(-)}\right)\theta}\left(3\left|iii\right\rangle
-\left|\bar{i}\bar{i}i\right\rangle-\left|i\bar{i}\bar{i}\right\rangle-\left|\bar{i}i\bar{i}\right\rangle\right)\nonumber\\
&&{\bf T}^{(3)} \left(\left|\bar{i}\bar{i}i\right\rangle+e^{{\bf
i}\frac{2\pi}3}\left|i\bar{i}\bar{i}\right\rangle+e^{{\bf
i}\frac{4\pi}3}\left|\bar{i}i\bar{i}\right\rangle\right)
=e^{{\bf i}\frac{2\pi}3}e^{\left(m_{ii}^{(+)}+2m_{ii}^{(-)}\right)\theta}\left(\left|\bar{i}\bar{i}i\right\rangle
+e^{{\bf i}\frac{2\pi}3}\left|i\bar{i}\bar{i}\right\rangle+e^{{\bf i}\frac{4\pi}3}\left|\bar{i}i\bar{i}\right\rangle\right)\nonumber\\
&&{\bf T}^{(3)} \left(\left|\bar{i}\bar{i}i\right\rangle+e^{{\bf
i}\frac{4\pi}3}\left|i\bar{i}\bar{i}\right\rangle+e^{{\bf
i}\frac{2\pi}3}\left|\bar{i}i\bar{i}\right\rangle\right) =e^{{\bf
i}\frac{4\pi}3}e^{\left(m_{ii}^{(+)}+2m_{ii}^{(-)}\right)\theta}\left(\left|\bar{i}\bar{i}i\right\rangle+e^{{\bf
i}\frac{4\pi}3}\left|i\bar{i}\bar{i}\right\rangle+e^{{\bf
i}\frac{2\pi}3}\left|\bar{i}i\bar{i}\right\rangle\right).\nonumber\end{eqnarray}
For $\overline{V_\lambda}$ one follows exactly similar steps and
obtains the same set of eigenvalues (with $i\rightleftharpoons
\bar{i}$ in the eigenstates). Combining all the solutions above we
obtain complete results for $r=3$ and all even $N$.

\paragraph{$\bullet$ $r=4$:} The key result (3.8) indicates the following construction (a direct generalization of (3.17)).
Define, for the even subset,
\begin{eqnarray}
&&V_{\left(b_1b_2b_3b_4\right)}=\left|b_1b_2b_3b_4\right\rangle+\epsilon_{12}\left|\overline{b_1}\overline{b_2}b_3b_4\right\rangle+
\epsilon_{23}\left|b_1\overline{b_2}\overline{b_3}b_4\right\rangle+\epsilon_{34}\left|b_1b_2\overline{b_3}\overline{b_4}\right\rangle+\\
&&\phantom{V_{b_1b_2b_3b_4}=}\epsilon_{41}\left|\overline{b_1}b_2b_3\overline{b_4}\right\rangle+
\epsilon_{12}\epsilon_{23}\left|\overline{b_1}b_2\overline{b_3}b_4\right\rangle+\epsilon_{12}\epsilon_{34}\left|\overline{b_1}\overline{b_2}\overline{b_3}\overline{b_4}\right\rangle
+\epsilon_{23}\epsilon_{34}\left|b_1\overline{b_2}b_3\overline{b_4}\right\rangle\nonumber
\end{eqnarray}
and similarly, implementing circular permutations as in (3.17),
$V_{\left(b_4b_1b_2b_3\right)},V_{\left(b_3b_4b_1b_2\right)},
V_{\left(b_2b_3b_4b_1\right)}$. Define
\begin{equation}
V_{\lambda}=V_{\left(b_1b_2b_3b_4\right)}+\lambda
V_{\left(b_4b_1b_2b_3\right)}+\lambda^2
V_{\left(b_3b_4b_1b_2\right)}+\lambda^3V_{\left(b_2b_3b_4b_1\right)},
\end{equation}
where
\begin{equation}
\lambda=\left(1,e^{{\bf i}\frac{2\pi}4} ,e^{{\bf
i}2\cdot\frac{2\pi}4},e^{{\bf i}3\cdot\frac{2\pi}4}\right)
=\left(1,{\bf i} ,-1,-{\bf i}\right).
\end{equation}
One obtains for distinct $\left(i_1,i_2,i_3,i_4\right)$
\begin{equation}
{\bf T}^{(4)} V_{\lambda}=\lambda
e^{\left(m_{i_1i_2}^{(\epsilon_{12})}+m_{i_2i_3}^{(\epsilon_{23})}+m_{i_3i_4}^{(\epsilon_{34})}+m_{i_4i_1}^{(\epsilon_{41})}
\right)\theta}V_{\lambda},
\end{equation}
where for non-zero results $\left(1+\epsilon_{12}\epsilon_{23}
\epsilon_{34}\epsilon_{41}\right)=2$. For the 8 possibilities
remaining for the $\epsilon$'s and the 4 values of $\lambda$ one
obtains the full set of 32 eigenstates and the corresponding
eigenvalues. For the complementary "odd" subspace formed by
\begin{eqnarray}
&&\left\{\left|b_1b_2b_3\overline{b_4}\right\rangle,\,\left|b_1b_2\overline{b_3}b_4\right\rangle,\,
\left|b_1\overline{b_2}b_3b_4\right\rangle,\,\left|\overline{b_1}b_2b_3b_4\right\rangle,\right.\nonumber\\
&&\left.\left|\overline{b_1}\overline{b_2}\overline{b_3}b_4\right\rangle,
\,\left|\overline{b_1}\overline{b_2}b_3\overline{b_4}\right\rangle,\,\left|\overline{b_1}b_2\overline{b_3}\overline{b_4}\right\rangle
,\,\left|b_1\overline{b_2}\overline{b_3}\overline{b_4}\right\rangle\right\}
\end{eqnarray}
and the sets related through circular permutations of the indices
one obtains, entirely in evident analogy of (3.25),
\begin{equation}
\overline{V_{\lambda}}=\overline{V_{\left(b_1b_2b_3b_4\right)}}+\lambda
\overline{V_{\left(b_4b_1b_2b_3\right)}}+\lambda^2
\overline{V_{\left(b_3b_4b_1b_2\right)}}+\lambda^3\overline{V_{\left(b_2b_3b_4b_1\right)}}.
\end{equation}
One obtains, for distinct indices, the same set of eigenvalues as
in (3.27). For some indices equal, but not all, the preceding
results can be carried over essentially. For all indices equal one
obtains the following spectrum of eigenvalues (grouping together
the even and the odd subspaces):
\begin{eqnarray}
&& \phantom{\hbox{[3-times]:}\,\,\,}e^{4m_{ii}^{(+)}\theta}\left(1,1\right)\nonumber\\
&&\hbox{[3-times]:}\,\,\,e^{2\left(m_{ii}^{(+)}+m_{ii}^{(-)}\right)\theta}\left(1,e^{{\bf
i}\frac{2\pi}4} ,e^{{\bf i}2\cdot\frac{2\pi}4},e^{{\bf i}3\cdot\frac{2\pi}4}\right)\nonumber\\
&&
\phantom{\hbox{[3-times]:}\,\,\,}e^{4m_{ii}^{(-)}\theta}\left(1,-1\right)
\end{eqnarray}
(In fact these are the same as in (A.13) of Ref. 7, since the
distinguishing index for odd dimensions signalled in (2.14) and
(2.15) is not involved.) We will not present here the 16
orthogonal eigenstates corresponding to (3.30). But a new feature,
as compared to (3.23), should be pointed out. For $r=3$,
considering together even and odd spaces, one has for identical
indices the eigenvalues
\begin{eqnarray}
&& \phantom{\hbox{[3-times]:}\,\,\,}e^{3m_{ii}^{(+)}\theta}\left(1,1\right)\nonumber\\
&&\hbox{[2-times]:}\,\,\,e^{\left(m_{ii}^{(+)}+2m_{ii}^{(-)}\right)\theta}\left(1,e^{{\bf
i}\frac{2\pi}3} ,e^{{\bf i}2\cdot\frac{2\pi}3}\right).
\end{eqnarray}
In fact, for each $r$ one obtains for $\epsilon=1$,
$e^{rm_{ii}^{(+)}\theta}\left(1,1\right)$ giving
\begin{equation}
\hbox{Tr}\left({\bf
T}^{(r)}\right)=2\sum_{i=1}^{n}e^{rm_{ii}^{(\epsilon)}\theta}.
\end{equation}
All other multiplets have each one zero sum. Concerning the zero
sum multiplets, comparing (3.31) and (3.30), we note that:

\begin{enumerate}
    \item For $r=3$ one has only triplets.
    \item For $r=4$ one has quadruplets and also a doublet
    $\left(1,-1\right)$. This last enters due to factorizability
    of $r$ ($4=2\times 2$).
\end{enumerate}
Such a feature is worth signalling since it generalizes.

When $r$ is a prime number the zero-sum multiplets appear only as
$r$-plets (with $r$-th roots of unity). How they consistently
cover the whole base space (along with the doublets giving the
trace as in (3.32)) has been amply discussed in our previous
papers \cite{R7,R10} under the heading "An encounter with a
theorem of Fermat". When $r$ has many factors the multiplicity of
sub-multiplets corresponding to each one is difficult to formulate
giving, say, an explicit general prescription.

\paragraph{$\bullet$ $r\geq 5$:} Another problem concerning
systematic, complete enumeration one encounters already for $r=5$.
Starting with (3.9) the generalizations of (3.24-27) are quite
evident. We will not present them here. But as compared to(3.30)
one has now for the 32-dimensional subspace (for
$i_1=i_2=\cdots=i_5$)
\begin{eqnarray}
&& \phantom{\hbox{[$n_3$-times]:}}\,\,\,e^{5m_{ii}^{(+)}\theta}\left(1,1\right)\nonumber\\
&&\hbox{[$n_2$-times]:}\,\,\,e^{\left(3m_{ii}^{(+)}+2m_{ii}^{(-)}\right)\theta}\left(1,e^{{\bf
i}\frac{2\pi}5} ,e^{{\bf i}2\cdot\frac{2\pi}5},e^{{\bf i}3\cdot\frac{2\pi}5},e^{{\bf i}4\cdot\frac{2\pi}5}\right)\nonumber\\
&&\hbox{[$n_4$-times]:}\,\,\,e^{\left(m_{ii}^{(+)}+4m_{ii}^{(-)}\right)\theta}\left(1,e^{{\bf
i}\frac{2\pi}5} ,e^{{\bf i}2\cdot\frac{2\pi}5},e^{{\bf
i}3\cdot\frac{2\pi}5}, e^{{\bf i}4\cdot\frac{2\pi}5}\right),
\end{eqnarray}
where $n_2+n_4=6$. At this stage it is not difficult to present
complete construction. But as $r$ increases, for prime numbers,
one has
\begin{eqnarray}
&& \phantom{\hbox{[$n_{r-1}$-times]:}}\,\,\,e^{rm_{ii}^{(+)}\theta}\left(1,1\right)\nonumber\\
&&\hbox{[$n_2$-times]:}\,\,\,\,\,\,\,e^{\left(\left(r-2\right)m_{ii}^{(+)}+2m_{ii}^{(-)}\right)\theta}\left(1,e^{{\bf
i}\frac{2\pi}r} ,\ldots,e^{{\bf i}\left(r-1\right)\cdot\frac{2\pi}r}\right)\nonumber\\
&&\hbox{[$n_4$-times]:}\,\,\,\,\,\,\,e^{\left(\left(r-4\right)m_{ii}^{(+)}+4m_{ii}^{(-)}\right)\theta}\left(1,e^{{\bf
i}\frac{2\pi}r} ,\ldots,e^{{\bf i}\left(r-1\right)\cdot\frac{2\pi}r}\right),\nonumber\\
&&\phantom{xxxx}\vdots\phantom{xxxxxxxxxxxxxxxxxxxxxxx}\vdots\nonumber\\
&&\hbox{[$n_{r-1}$-times]:}\,\,\,e^{\left(m_{ii}^{(+)}+\left(r-1\right)m_{ii}^{(-)}\right)\theta}\left(1,e^{{\bf
i}\frac{2\pi}r} ,\ldots,e^{{\bf
i}\left(r-1\right)\cdot\frac{2\pi}r}\right)
\end{eqnarray}
the sum of the multiplicities of such $r$-plets satisfying
\begin{equation}
n_2+n_4+\cdots +n_{r-1}=\frac 1r\left(2^r-2\right).
\end{equation}
(That the right hand side is an integer for $r$ a prime number is
guaranteed by a theorem of Fermat, as has been discussed before.)
A general explicit prescription for the sequence
$\left(n_2,n_4,\ldots,n_{r-1}\right)$ is beyond the scope of this
paper.

The situation is as above for prime numbers $r$. For factorizable
$r$ the presence, in addition, of submultiplets has already been
pointed out. For $r=4$, only such submultiplets was
$\left(1,-1\right)$. For $r=6$, one can have $\left(1,-1\right)$,
$\left(1,e^{{\bf i}\frac{2\pi}3},e^{{\bf
i}2\cdot\frac{2\pi}3}\right)$ corresponding to the factors
$6=2\times 3$ respectively. For $r=30=2\times 3\times 5$ one can
have also submultiplets corresponding to 5. And so on.

But apart from the above mentioned limitations concerning
multiplicities of zero-sum ("roots of unity") multiplets and
submultiplets we can claim to have elucidated the spectrum of
eigenvalues for all $\left(N,r\right)$. For the generic case with
distinct indices for any $\left(N,r\right)$, one has eigenvalues
\begin{equation}
\lambda_ke^{\left(m_{i_1i_2}^{\left(\epsilon_{12}\right)}+m_{i_2i_3}^{\left(\epsilon_{23}\right)}+\cdots+
m_{i_{r-1}i_r}^{\left(\epsilon_{r-1,r}\right)}+
m_{i_1i_r}^{\left(\epsilon_{1r}\right)}\right)\theta},
\end{equation}
where $\left(1+\epsilon_{12}\epsilon_{23}\cdots
\epsilon_{r-1,r}\epsilon_{1,r}\right)=2$ and $\lambda_k=e^{{\bf
i}k\cdot\frac{2\pi}r}$ $\left(k=0,1,2,\ldots,r-1\right)$. Here
each $m_{i_ki_l}^{\left(\epsilon_{kl}\right)}$ is a free
parameter. We solve only sets of linear equations with quite
simple constant coefficients to obtain the eigenstates and the
eigenvalues.

\paragraph{\bf B. Odd dimensions:} We refer back to (2.14-16). When
the special index $n$ ($n=\overline{n}$) for $N=2n-1$ is not
present in the basis state $\left|b_1b_2\ldots b_r\right\rangle$
and hence in the subspace it generates via cyclic permutations as
in the forgoing examples, the foregoing constructions can be taken
over wholesale. When $n$ is present the modifications are not
difficult to take into account. The trace is now
$\hbox{Tr}\left({\bf
T}^{(r)}\right)=2\sum_{i=1}^{n-1}e^{rm_{ii}^{(\epsilon)}\theta}+1$,
since the central state $\left|n\right\rangle\otimes \left|n\right
\rangle\otimes\ldots \left|n\right\rangle\otimes\equiv \left|nn
\ldots n\right\rangle$ contributes with our normalization 1 to the
trace. Various aspects have been studied in considerable detail in
our previous papers \cite{R6,R7} on odd $N$. Here we just refer to
them.

\section{Spin Chain Hamiltonians}
\setcounter{equation}{0}

Spin chains corresponding to our braid matrices have already been
studied in our previous papers \cite{R7,R9}. Here we formulate a
unified approach for all dimensions $\left(N=2n-1,2n\right)$.

The basic formula (see sources cited in Refs. 7, 8) is
\begin{equation}
H=\sum_{k=1}^rI\otimes\cdots\otimes I\otimes \dot{\hat{R}}_{k,k+1}
\left(0\right)\otimes I\otimes\cdots\otimes I,
\end{equation}
where for circular boundary conditions,  for $k+1=r+1\approx 1$.
For even $N$ ($N=2n$), from (2.4),
\begin{equation}
\dot{\widehat{R}}\left(0\right)=\left.\frac{d}{d\theta}
\hat{R}\left(\theta\right)\right|_{\theta=0}=\sum_{\epsilon,i,j}
m_{ij}^{(\epsilon)}\left(P_{ij}^{(\epsilon)}+P_{i\overline{j}}^{(\epsilon)}\right)
\end{equation}
the projectors being given by (2.1) and the remark below (2.1).
For $N$ odd ($N=2n-1$; $n=2,3,\ldots$) we introduce a modified
overall normalization factor to start with. Such a factor is
trivial concerning the braid equation, but not for the Hamiltonian
(since a derivative is involved) if the factor is
$\theta$-dependent. Multiply (2.16) by $e^{m\theta}$ and redefine
\begin{equation}
\left(\left(m+m_{ij}^{(\epsilon)}\right),\left(m+m_{ni}^{(\epsilon)}\right),\left(m+m_{in}^{(\epsilon)}\right)\right)
\longrightarrow\left(m_{ij}^{(\epsilon)},m_{ni}^{(\epsilon)},m_{in}^{(\epsilon)}\right)
\end{equation}
since $\left(m_{ij}^{(\epsilon)},\ldots\right)$ are arbitrary to
start with. Now for odd $N$ (with the ranges of
$\left(i,j\right)$, $\left(\bar{i},\bar{j}\right)$ of (2.16))
\begin{equation}
\dot{\widehat{R}}\left(0\right)=mP_{nn}+\sum_{\epsilon,i}
\left(m_{ni}^{(\epsilon)}P_{ni}^{(\epsilon)}+m_{in}^{(\epsilon)}P_{in}^{(\epsilon)}\right)+\sum_{\epsilon,i,j}
m_{ij}^{(\epsilon)}\left(P_{ij}^{(\epsilon)}+P_{i\overline{j}}^{(\epsilon)}\right).
\end{equation}
This extends the $m=0$ case by including the presence of $P_{nn}$.

In $H$, $\dot{\widehat{R}}\left(0\right)_{k,k+1}$ acts on the
basis $\left|V\right\rangle_{(k)}\otimes\left|V\right
\rangle_{(k+1)}$. One can denote, using standard ordering of spin
components for each $k$,
\begin{equation}
\left|V\right\rangle_{(k)}=\left|\begin{matrix}
  \left|n-1/2\right\rangle_{k} \\
  \vdots \\
  \left|1/2\right\rangle_{k} \\
  \left|-1/2\right\rangle_{k} \\
  \vdots \\
  \left|-n+1/2\right\rangle_{k} \\
\end{matrix}\right\rangle,\qquad \left|V\right\rangle_{(k)}=\left|\begin{matrix}
  \left|n-1\right\rangle_{k} \\
  \vdots \\
  \left|1\right\rangle_{k} \\
\left|0\right\rangle_{k} \\
  \left|-1\right\rangle_{k} \\
  \vdots \\
  \left|-n+1\right\rangle_{k} \\
\end{matrix}\right\rangle
\end{equation}
for $N=2n,\,2n-1$ respectively. Without being restricted to spin
one can consider more generally any system with $N$ orthogonal
states. We will continue however to use the terminology of spin.
At each site one can consider a superposition at each level such
as $\left(\sum_lc_{il}^{(k)}\left|l\right\rangle_k\right)$. To
keep the notation tractable we will just denote
\begin{equation}
\left|V\right\rangle_{(k)}=\left|\begin{matrix}
  \left|1\right\rangle_{k} \\
  \left|2\right\rangle_{k} \\
  \vdots \\
  \left|\bar{2}\right\rangle_{k} \\
  \left|\bar{1}\right\rangle_{k} \\
\end{matrix}\right\rangle
\end{equation}
and keep possible significances of $\left|i\right\rangle_k$ in
mind. Each index $k$ is acted upon twice by $H$, namely by
$\dot{\hat{R}}_{k-1,k} \left(0\right)$, $\dot{\hat{R}}_{k,k+1}
\left(0\right)$ and, for closed chains, the index 1 is also thus
involved in $\dot{\hat{R}}_{12} \left(0\right)$ and
$\dot{\hat{R}}_{r1} \left(0\right)$. Some simple examples are:
\begin{description}
    \item[(i)] $N=2$:
\begin{equation}
\dot{\hat{R}}\left(0\right)=\begin{vmatrix}
  \hat{a}_+ & 0 & 0 & \hat{a}_- \\
  0 & \hat{a}_+ & \hat{a}_- & 0 \\
  0 & \hat{a}_- & \hat{a}_+ & 0 \\
  \hat{a}_- & 0 & 0 & \hat{a}_+ \\
\end{vmatrix},
\end{equation}
where $\hat{a}_{\pm}=\frac 12\left(m_{11}^{(+)}\pm
m_{11}^{(-)}\right)$ and hence in evident notations
\begin{equation}
\dot{\hat{R}}\left(0\right)\left|V\right\rangle_{(k)}\otimes\left|V\right
\rangle_{(k+1)}=\left|\begin{matrix}
  \hat{a}_+\left|11\right\rangle +\hat{a}_-\left|\bar{1}\bar{1}\right\rangle \\
  \hat{a}_+\left|1\bar{1}\right\rangle +\hat{a}_-\left|\bar{1}1\right\rangle \\
  \hat{a}_-\left|1\bar{1}\right\rangle +\hat{a}_+\left|\bar{1}1\right\rangle \\
  \hat{a}_-\left|11\right\rangle +\hat{a}_+\left|\bar{1}\bar{1}\right\rangle \\
\end{matrix}\right\rangle_{\left(k,k+1\right)},
\end{equation}

    \item[(ii)] $N=4$: In the notation of sec. 7 of Ref. 9
    \begin{equation}
\hat{R}\left(\theta\right)=\begin{vmatrix}
  \hat{D}_{11} & 0 & 0 & \hat{A}_{1\bar{1}} \\
  0 & \hat{D}_{22} & \hat{A}_{2\bar{2}} & 0 \\
  0 & \hat{A}_{\bar{2}2} & \hat{D}_{\bar{2}\bar{2}} & 0 \\
  \hat{A}_{\bar{1}1} & 0 & 0 & \hat{D}_{\bar{1}\bar{1}} \\
\end{vmatrix},
\end{equation}
where
\begin{eqnarray}
&&\hat{D}_{11}=\hat{D}_{\bar{1}\bar{1}}=\left(\begin{array}{cccc}
  \hat{a}_+ & 0 & 0 & 0 \\
  0 & \hat{b}_+ & 0 & 0 \\
  0 & 0 & \hat{b}_+ & 0 \\
  0 & 0 & 0 & \hat{a}_+ \\
\end{array}\right),\qquad
\hat{D}_{22}=\hat{D}_{\bar{2}\bar{2}}=\left(\begin{array}{cccc}
  \hat{c}_+ & 0 & 0 & 0 \\
  0 & \hat{d}_+ & 0 & 0 \\
  0 & 0 & \hat{d}_+ & 0 \\
  0 & 0 & 0 & \hat{c}_+ \\
\end{array}\right),\nonumber\\
&&\hat{A}_{1\bar{1}}=\hat{A}_{\bar{1}1}=\left(\begin{array}{cccc}
  0 & 0 & 0 & \hat{a}_- \\
  0 & 0 & \hat{b}_- & 0 \\
  0 & \hat{b}_- & 0 & 0 \\
  \hat{a}_- & 0 & 0 & 0 \\
\end{array}\right),\qquad
\hat{A}_{2\bar{2}}=\hat{A}_{\bar{2}2}=\left(\begin{array}{cccc}
  0 & 0 & 0 & \hat{c}_- \\
  0 & 0 & \hat{d}_- & 0 \\
  0 & \hat{d}_- & 0 & 0 \\
  \hat{c}_- & 0 & 0 & 0 \\
\end{array}\right),\nonumber\\
&&\hat{a}_{\pm}=\frac 12\left(m_{11}^{(+)}\pm
m_{11}^{(-)}\right),\qquad \hat{b}_{\pm}=\frac
12\left(m_{12}^{(+)}\pm m_{12}^{(-)}\right),\nonumber\\
&&\hat{c}_{\pm}=\frac 12\left(m_{21}^{(+)}\pm
m_{21}^{(-)}\right),\qquad \hat{d}_{\pm}=\frac
12\left(m_{22}^{(+)}\pm m_{22}^{(-)}\right)\end{eqnarray} Hence,
\begin{equation}
\dot{\hat{R}}\left(0\right)\left|V\right\rangle_{(k)}\otimes\left|V\right
\rangle_{(k+1)}=\left|\begin{matrix}
  \left|1\right\rangle\otimes\left(\hat{D}_{11}\left|V\right\rangle\right)+
  \left|\bar{1}\right\rangle\otimes\left(\hat{A}_{1\bar{1}}\left|V\right\rangle\right) \\
  \left|2\right\rangle\otimes\left(\hat{D}_{22}\left|V\right\rangle\right)+
  \left|\bar{2}\right\rangle\otimes\left(\hat{A}_{2\bar{2}}\left|V\right\rangle\right) \\
  \left|2\right\rangle\otimes\left(\hat{A}_{2\bar{2}}\left|V\right\rangle\right)+
  \left|\bar{2}\right\rangle\otimes\left(\hat{D}_{22}\left|V\right\rangle\right) \\
  \left|1\right\rangle\otimes\left(\hat{A}_{1\bar{1}}\left|V\right\rangle\right)+
  \left|\bar{1}\right\rangle\otimes\left(\hat{D}_{11}\left|V\right\rangle\right) \\
\end{matrix}\right\rangle_{\left(k,k+1\right)},
\end{equation}
where for each $k$, $\left|V\right\rangle=\left|\begin{matrix}
  \left|1\right\rangle\\
  \left|2\right\rangle\\
  \left|\bar{2}\right\rangle\\
  \left|\bar{1}\right\rangle\\
\end{matrix}\right\rangle$.

\item[(iii)] $N=3$: Adapting the results of Sec. 1 of ref. 7 and
Sec. 11 of Ref. 8 to notations analogous to the cases above, one
can write
\begin{equation}
\hat{R}\left(\theta\right)=\begin{vmatrix}
  D & 0 & A \\
  0 & C & 0 \\
  A & 0 & D \\
  \end{vmatrix},
\end{equation}
where
\begin{equation}
D=\left(\begin{array}{ccc}
  \hat{a}_+ & 0 & 0 \\
  0 & \hat{b}_+ & 0 \\
  0 & 0 & \hat{a}_+ \\
  \end{array}\right),\qquad
A=\left(\begin{array}{ccc}
  0 & 0 & \hat{a}_- \\
  0 & \hat{b}_- & 0 \\
  \hat{a}_- & 0 & 0 \\
  \end{array}\right),\qquad
C=\left(\begin{array}{ccc}
  \hat{c}_+ & 0 & \hat{c}_- \\
  0 & m & 0 \\
  \hat{c}_- & 0 & \hat{c}_+ \\
\end{array}\right)
\end{equation}
the central element $m$ corresponding to $mP_{nn}$ of (4.4) and
\begin{equation}
\hat{a}_{\pm}=\frac 12\left(m_{11}^{(+)}\pm
m_{11}^{(-)}\right),\qquad \hat{b}_{\pm}=\frac
12\left(m_{12}^{(+)}\pm m_{12}^{(-)}\right),\qquad
\hat{c}_{\pm}=\frac 12\left(m_{21}^{(+)}\pm
m_{21}^{(-)}\right),\end{equation}
\end{description}
Denote the basis factors for each $k$ as,
$\left|V\right\rangle_{(k)}=\left|\begin{matrix}
\left|1\right\rangle\\\left|2\right\rangle\\\left|\bar{1}\right\rangle\\
\end{matrix}\right\rangle_{(k)}$ and
\begin{equation}
\left|V\right\rangle_{(k)}\otimes\left|V\right
\rangle_{(k+1)}=\left|\begin{matrix}
  \left|1\right\rangle \otimes\left|V\right\rangle \\
  \left|2\right\rangle \otimes\left|V\right\rangle \\\left|\bar{1}\right\rangle \otimes\left|V\right\rangle \\
\end{matrix}\right\rangle_{\left(k,k+1\right)}.
\end{equation}
One obtains
\begin{equation}
\dot{\hat{R}}\left(0\right)\left|V\right\rangle_{(k)}\otimes\left|V\right
\rangle_{(k+1)}=\begin{vmatrix}
  D & 0 & A \\
  0 & C & 0 \\
  A & 0 & D \\
  \end{vmatrix}\left|\begin{matrix}
  \left|1\right\rangle \otimes\left|V\right\rangle \\
  \left|2\right\rangle \otimes\left|V\right\rangle \\\left|\bar{1}\right\rangle \otimes\left|V\right\rangle \\
\end{matrix}\right\rangle_{\left(k,k+1\right)}
=\left|\begin{matrix}
  \left|1\right\rangle \otimes\left(D\left|V\right\rangle\right)+\left|\bar{1}\right\rangle \otimes\left(A\left|V\right\rangle\right) \\
\left|2\right\rangle \otimes\left(C\left|V\right\rangle\right) \\
\left|1\right\rangle \otimes\left(A\left|V\right\rangle\right)+\left|\bar{1}\right\rangle \otimes\left(D\left|V\right\rangle\right) \\
\end{matrix}\right\rangle_{\left(k,k+1\right)}.
\end{equation}

For $N=5,6,\ldots$ one can generate such results systematically.
They, considering the action of all the terms of (4.1), furnish
the transition matrix elements and expectation values for possible
states of the chain and permit a study of correlations.

One can also consider higher order conserved quantities (Sec. 1.5
of Ref. 10) given by
\begin{equation}
H_l=\left.\frac{d^l}{d\theta^l} \log{\bf
T}^{(r)}\left(\theta\right)\right|_{\theta=0}.
\end{equation}
For $l=2$, (4.1) is generalized by the appearance of factors of
the type (apart from non-overlapping derivatives)
\begin{equation}
\dot{\hat{R}}\left(0\right)_{k-1,k}\otimes
\dot{\hat{R}}\left(0\right)_{k,k+1},\qquad
\ddot{\hat{R}}\left(0\right)_{k,k+1}.
\end{equation}
For $l>2$ this generalizes in an evident fashion. A study of spin
chains for the "exotic" S\O3 can be found in Ref. 11. It is
interesting to compare it with the $4\times 4$ (for $N=2$) case
briefly presented above by exploring the latter case in comparable
detail.

One may note that for our class of braid matrices, for even $N$,
as compared to (4.2),
\begin{equation}
\left.\frac{d^l}{d\theta^l}
\hat{R}\left(\theta\right)\right|_{\theta=0}=\sum_{\epsilon}\sum_{i,j}
\left(m_{ij}^{(\epsilon)}\right)^l\left(P_{ij}^{(\epsilon)}+P_{i\overline{j}}^{(\epsilon)}\right)
\end{equation}
and for odd $N$, as compared to (4.4)
\begin{equation}
\left.\frac{d^l}{d\theta^l}
\hat{R}\left(\theta\right)\right|_{\theta=0}=m^lP_{nn}+\sum_{\epsilon,i}
\left(\left(m_{ni}^{(\epsilon)}\right)^lP_{ni}^{(\epsilon)}+\left(m_{in}^{(\epsilon)}\right)^lP_{in}^{(\epsilon)}\right)+\sum_{\epsilon,i,j}
\left(m_{ij}^{(\epsilon)}\right)^l\left(P_{ij}^{(\epsilon)}+P_{i\overline{j}}^{(\epsilon)}\right).
\end{equation}

\section{Potentials for factorizable $S$-matrices}
\setcounter{equation}{0}

Such potentials can be obtained as inverse Cayley transforms of
Yang-Baxter matrices of appropriate dimensions (see for example,
sec. 3 of Ref. 12 and sec. 1 of Ref. 2). Starting with the
Yang-Baxter matrix $R\left(\theta\right)={\bf P}\hat{R}
\left(\theta\right)$ the required potential ${\bf
V}\left(\theta\right)$ is given by
\begin{equation}
-{\bf i}{\bf V}\left(\theta\right)=
\left(R\left(\theta\right)-\lambda\left(\theta\right)I\right)^{-1}\left(R\left(\theta\right)
+\lambda\left(\theta\right)I\right).
\end{equation}
We have emphasized in our previous studies (sec. 5 of Ref. 7, sec.
8 of Ref. 8) that $\lambda\left(\theta\right)$ cannot be
arbitrary, a set of values must be excluded for the inverse (5.1)
to be well-defined. Here we generalize our previous results to all
$N$. As compared to other well-known studies \cite{R14} of fields
corresponding to factorizable scatterings here our construction
starts with the braid matrices (2.4) and (2.16).

Define
\begin{equation}
X\left(\theta\right)=\left(R\left(\theta\right)-\lambda\left(\theta\right)I\right)^{-1}
\end{equation}
when
\begin{equation}
-{\bf i}{\bf V}\left(\theta\right)=I+2\lambda\left(\theta\right)
X\left(\theta\right).
\end{equation}
We give below immediately the general solution and then explain
the notations more precisely. The solution was obtained via formal
series expansions. But the final closed form can be verified
directly. The solution is,
\begin{eqnarray}
&&X\left(\theta\right)=-\frac 12 \sum_{\epsilon=\pm}
\sum_{a,b=1}^N\displaystyle
\frac{1}{\lambda^2\left(\theta\right)-e^{\left(m_{ab}^{(\epsilon)}+m_{ba}^{(\epsilon)}\right)\theta}}
\left\{\lambda\left(\theta\right)\left[\left(aa\right)\otimes
\left(bb\right)+ \epsilon \left(a\overline{a}\right)\otimes
\left(b\overline{b}\right)\right]+\right.\nonumber\\
&&\phantom{X\left(\theta\right)=}\left.e^{m_{ba}^{(\epsilon)}\theta}
\left[\left(ab\right)\otimes \left(ba\right)+ \epsilon
\left(a\overline{b}\right)\otimes
\left(b\overline{a}\right)\right] \right\}
\end{eqnarray}
where
\begin{equation}
\lambda\left(\theta\right)\neq \pm e^{\frac
12\left(m_{ab}^{(\epsilon)}+m_{ba}^{(\epsilon)}\right)\theta}.
\end{equation}
For $N=2n$, (2.9) is implicit in this result. For $N=2n-1$, when
$\overline{n}=n$ the conventions (2.16) (or (4.4)) are to be
implemented when $a$ or $b$ or both $\left(a,b\right)$ are $n$.
From (5.3) and (5.4),
\begin{eqnarray}
&&{\bf V}\left(\theta\right)=-\frac {\bf i}2 \sum_{\epsilon=\pm}
\sum_{a,b=1}^N\displaystyle
\left\{\frac{\lambda^2\left(\theta\right)+e^{\left(m_{ab}^{(\epsilon)}+m_{ba}^{(\epsilon)}\right)\theta}}
{\lambda^2\left(\theta\right)-e^{\left[m_{ab}^{(\epsilon)}+m_{ba}^{(\epsilon)}\right]\theta}}
\left[\left(aa\right)\otimes \left(bb\right)+ \epsilon
\left(a\overline{a}\right)\otimes\left(b\overline{b}\right)\right]+\right.\nonumber\\
&&\phantom{X\left(\theta\right)=}\left.\frac{e^{m_{ab}^{(\epsilon)}\theta}}
{\lambda^2\left(\theta\right)-e^{\left(m_{ab}^{(\epsilon)}+m_{ba}^{(\epsilon)}\right)\theta}}
\left[\left(ba\right)\otimes \left(ab\right)+ \epsilon
\left(b\overline{a}\right)\otimes
\left(a\overline{b}\right)\right]
\right\}\nonumber\\
&&\phantom{X\left(\theta\right)}\equiv  \sum_{ab,cd}{\bf
V}\left(\theta\right)_{\left(ab,cd\right)}\left(ab\right)\otimes\left(cd\right).
\end{eqnarray}
Now one can write down the Lagrangians for scalar and spinor
fields. For the spinor case, for example,
\begin{equation}
{\cal L}=\int dx\left[{\bf
i}\overline{\psi}_a\gamma_\nu\partial_\nu\psi_a
-\textsf{g}\left(\overline{\psi}_a\gamma_\nu\psi_c\right) {\bf
V}_{ab,cd} \left(\overline{\psi}_b\gamma_\nu\psi_d\right)\right].
\end{equation}
The simpler scalar case can be written analogously. The scalar
Lagrangian has an interaction term of the form
$\left(\overline{\phi}_a\overline{\phi}_c\right){\bf V}_{ab,cd}
\left(\phi_b\phi_d\right)$. For $N=2n$, one obtains the
nonvanishing elements of $\textbf{V}$ as
\begin{eqnarray}
&&{\bf
V}_{bb,dd}=-\frac{\textbf{i}}2\sum_{\epsilon}\frac{\lambda^2\left(\theta\right)+
e^{\left(m_{bd}^{(\epsilon)}+m_{db}^{(\epsilon)}\right)\theta}}
{\lambda^2\left(\theta\right)-e^{\left(m_{bd}^{(\epsilon)}+m_{db}^{(\epsilon)}\right)\theta}},\nonumber\\
&&{\bf
V}_{\overline{b}b,\overline{d}d}=-\frac{\textbf{i}}2\sum_{\epsilon}\epsilon\frac{\lambda^2\left(\theta\right)+
e^{\left(m_{bd}^{(\epsilon)}+m_{db}^{(\epsilon)}\right)\theta}}
{\lambda^2\left(\theta\right)-e^{\left(m_{bd}^{(\epsilon)}+m_{db}^{(\epsilon)}\right)\theta}},\nonumber\\
&&{\bf V}_{db,bd}=-\frac{\textbf{i}}2\sum_{\epsilon}\frac{
e^{m_{bd}^{(\epsilon)}\theta}}
{\lambda^2\left(\theta\right)-e^{\left(m_{bd}^{(\epsilon)}+m_{db}^{(\epsilon)}\right)\theta}},\nonumber\\
&&{\bf
V}_{\overline{d}b,\overline{b}d}=-\frac{\textbf{i}}2\sum_{\epsilon}\epsilon\frac{
e^{m_{bd}^{(\epsilon)}\theta}}
{\lambda^2\left(\theta\right)-e^{\left(m_{bd}^{(\epsilon)}+m_{db}^{(\epsilon)}\right)\theta}},
\end{eqnarray}
where $a,b\in\left\{1,\ldots,N\right\}$ and $\left((a,b\right)
\neq\left(n,n\right)$ if $N=2n-1$. When $N$ is odd ($N=2n-1$), one
has for the special index $n$,
\begin{equation}
{\bf
V}_{nn,nn}=-\textbf{i}\left(\frac{\lambda^2\left(\theta\right)+ 2}
{\lambda^2\left(\theta\right)-1}\right),
\end{equation}
where $m_{nn}^{(\epsilon)}$ are taken to be zero. For our models,
the scattering process can be schematically, presented as (see
figure)
\begin{figure}[ht]
\centerline{\includegraphics[height=7cm]{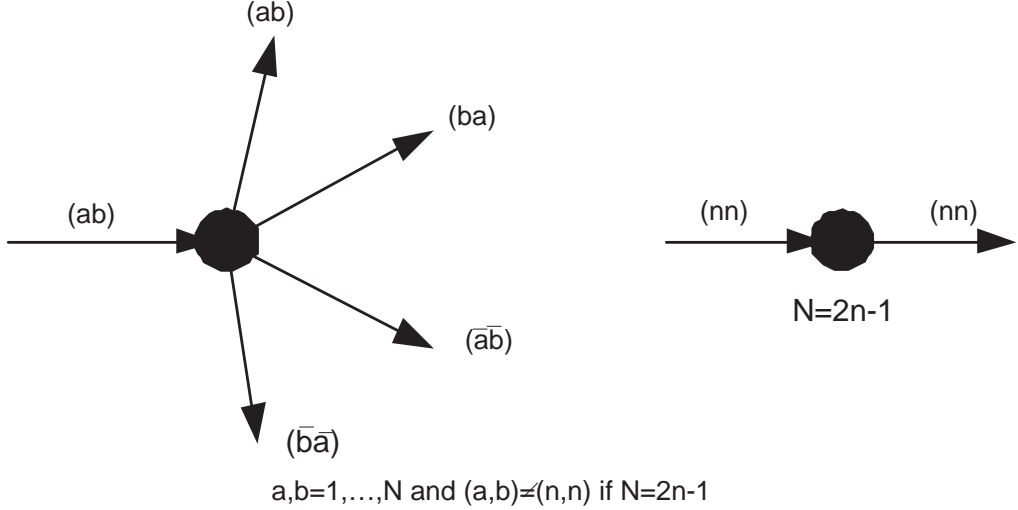}}
\caption{Scattering process}
\end{figure}

\section{Remarks}
\setcounter{equation}{0}

We have obtained explicit eigenvalues of the transfer matrix ${\bf
T}^{(r)}$ corresponding to our class of $N^2\times N^2$ braid
matrices for all $\left(r,N\right)$. Starting with our nested
sequence of projectors we obtained a very specific structure of
${\bf T}^{(r)}$. Exploiting this structure fully we obtained the
eigenvalues and eigenstates. The zero-sum multiplets parametrized
by sets of roots of unity arose from the circular permutations of
the indices in the tensor products of the basis states under the
action of ${\bf T}^{(r)}$. The same structure led to systematic
constructions of multiparameter Hamiltonians of spin chains
related to our class of braid matrices for all $N$. Finally we
constructed the inverse Cayley transformation for the general case
(any $N$) giving potentials compatible with factorizability of
$S$-matrices. The contents of such matrices (generalizing sec. 5
of Ref. 7) will be further studied elswhere. In the treatment of
each aspect we emphasized the role of a remarkable feature of our
formalism -- the presence of free parameters whose number increase
as  $N^2$ with $N$. A single class of constraints ($\theta\geq 0$,
$m_{ab}^{(+)}>m_{ab}^{(-)}$, $a,b\in\left\{1,\ldots,N\right\}$)
assures non-negative Boltzmann weights in the statistical models.
But our constructions also furnish directly some important
properties of such models. Thus $\hbox{Tr}\left({\bf
T}^{(r)}\right) = 2\sum_{i=1}^ne^{rm_{ii}^{(+)}\theta}$ ($N=2n$)
and $\hbox{Tr}\left({\bf T}^{(r)}\right) =
2\sum_{i=1}^{n-1}e^{rm_{ii}^{(+)}\theta}+1$ ($N=2n-1$) (for the
normalization (2.16)). The eigenvalues can be ordered in magnitude
by choosing the order of values of the parameters. Thus choosing
$m_{11}^{(+)}>m_{22}^{(+)}>\cdots$ the largest eigenvalue of
$\textbf{T}^{(r)}$ is for $\theta>0$, say
$e^{rm_{11}^{(+)}\theta}$ the next largest is
$e^{rm_{22}^{(+)}\theta}$ and so on. We intend to study elsewhere
the properties of our models more thoroughly.

In previous papers \cite{R9,ETAN} we pointed out that for purely
imaginary parameters ($\textbf{i}m_{ab}^{(\pm)}$ with
$m_{ab}^{(\pm)}$ real) our braid matrices are all unitary,
providing an entire class with free parameters for all $N$. Here
we have not repeated this discussion. But the fact that they
generate parametrized entangled states is indeed of interest.

Unitary matrices provide valid transformations of a basis (of
corresponding dimension) of quantum states. One may ask the
following question: If such a matrix, apart from being unitary,
also satisfies the braid equation what consequence might be
implied? Link between quantum and topological entanglements have
been discussed by several authors \cite{R15,R16} cited in our
previous papers \cite{R9,ETAN}. For the braid property to be
relevant a triple tensor product $\left(V\otimes V\otimes
V\right)$ of basis space is essential. We hope to explore
elsewhere our multiparameter unitary matrices in such a context.

\vskip 1cm

\paragraph{Note added:} Prof. J.H.H. Perk has kindly pointed out
that a class of multi-parameter generalization of the 6-vertex
model is provided by $sl\left(m|n\right)$ ones \cite{R17,R18}. As
Perk-Schultz models they have been studied by many authors and
have led to various important applications. (Relevant references
can be easily found via \textsc{arXiv}.) In this class the source of
parameters are multi-component rapidities. The study of
eigenvectors and eigenvalues were pioneered in Ref. 19; the most
recent follow-up being Ref. 20, where other references can be found.
Such studies may be compared to our systematic explicit constructions
for all $\left(r,N\right)$.

\vskip 0.5cm

\noindent{\bf Acknowledgments:} {\em One of us (BA) wants to thank
Pierre Collet and Paul Sorba for precious help. This work is
supported by the CNRS/DPGRF program number: 19841 "Nouvelles
solutions de Yang-Baxter et mod\`eles statistiques associ\'es;
M\'ethodes de quantifiations (super)jordaniennes des
(super)alg\`ebres; Cosmologie des Branes". }

\end{document}